\documentclass[11pt]{article}
\usepackage{amsmath,amssymb,theorem,amstext,amsgen,amsbsy,amsopn,amsfonts,graphicx,cases}
\usepackage{graphicx}
\usepackage{subfigure}
\usepackage{psfrag}
\usepackage{color}
\usepackage{enumerate}
\usepackage{epstopdf}
\usepackage{array}

\usepackage{latexsym}
\newtheorem{thm}{Theorem}[section]

\newtheorem{prop}[thm]{Proposition}
\newtheorem{ob}[thm]{Observation}
\theorembodyfont{\rmfamily}

\def\pf{\bigskip\noindent {\bf Proof.}}

\def\dfn#1{{\sl #1}}

\def\less{\setminus}

\def\pf{\bigskip\noindent {\emph{Proof.}}}

\def\qed{ \hfill\vrule height3pt width6pt depth2pt}

\def\pf{\bigskip\noindent {\bf Proof.  }}

\title{Planar  anti-Ramsey numbers of matchings}

\author{Gang Chen$^1$,  Yongxin Lan$^2$ and Zi-Xia Song$^3$\thanks{Corresponding Author.  Email address: Zixia.Song@ucf.edu}\\
$^1$School of mathematics and statistics\\
Ningxia University, China\\
 $^2$Center for Combinatorics and LPMC\\
Nankai University, Tianjin,  300071, China\\
$^3$Department  of Mathematics\\
 University of Central Florida\\
 Orlando, FL 32816, USA
}
\begin{document}
\maketitle
\begin{abstract}
 Given a positive integer $n$ and  a planar graph  $H$, let  $\mathcal{T}_n(H)$  be the family of  all plane triangulations $T$ on  $n$ vertices such that $T$  contains  a subgraph isomorphic to  $H$.  The \dfn{planar anti-Ramsey number of $H$}, denoted  $ar_{_\mathcal{P}}(n, H)$,  is the  maximum number  of colors in an edge-coloring of  a plane triangulation $T\in \mathcal{T}_n(H)$ such that $T$ contains   no    rainbow copy of $H$.   In this paper we study planar anti-Ramsey numbers of  matchings. For all $t\ge1$, let $M_t$ denote a matching of size $t$.    We  prove that for all  $t\ge6$ and  $n\ge 3t-6$,  $2n+3t-15\le ar_{_{\mathcal{P}}}(n,    {M}_t)\le  2n+4t-14$, which  significantly  improves the existing  lower and upper bounds for     $ar_{_\mathcal{P}}(n, M_t)$.  It seems that for each $t\ge6$, the lower bound we obtained   is the exact value  of $ar_{_{\mathcal{P}}}(n,    {M}_t)$ for sufficiently large $n$. This is indeed the case for $M_6$. We  prove that  $ar_{_\mathcal{P}}(n, M_6)=2n+3$ for all $n\ge30$.    \end{abstract}



 \baselineskip 15pt
\section{Introduction}
All graphs considered in this paper are finite and simple.  For a graph $G$ we use   $|G|$ and    $e(G)$  to  the  number
of vertices and    number of edges  of $G$, respectively.
For a vertex $x \in V(G)$, we will use $N_G(x)$ to denote the set of vertices in $G$ which are adjacent to $x$.
We define   $d_G(x) = |N_G(x)|$.
For any $A\subseteq V(G)$, the subgraph of $G$ induced by $A$, denoted $G[A]$, is the graph with vertex set $A$ and edge set $\{xy \in E(G) : x, y \in A\}$.   We denote   $G \less A$ the subgraph of $G$ induced on $V(G) \less A$.
If $A = \{a\}$, we simply write   $G \less a$. For disjoint subsets $A, B$ of $V(G)$,  we use $e_G(A, B)$ to denote the number of edges in $G$ with one end in $A$ and the other in $B$. Since every planar bipartite graph on $n\ge3$ vertices has at most $2n-4$ edges, we will frequently use the fact that $e_G(A, B)\le 2(|A|+|B|)-4$ when $G$ is planar and $|A\cup B|\ge3$.
Given two isomorphic graphs $G$ and $H$, we may (with a slight but common abuse of notation) write $G = H$. 
For any positive integer $k$, let  $[k]:=\{1,2, \ldots, k\}$. We use the convention that ``$A:=$'' means that $A$ is defined to be the right-hand side of the relation.\medskip

Motivated by anti-Ramsey numbers introduced by Erd\H{o}s, Simonovits  and S\'os~\cite{ESS} in 1975,  we study the anti-Ramsey problem when host graphs are plane triangulations.  
A subgraph of an edge-colored graph is \dfn{rainbow} if all of its edges have different colors.  Given    a planar graph  $H$ and a positive integer $n\ge|H|$, let  $\mathcal{T}_n(H)$  be the family of  all plane triangulations $T$ on  $n$ vertices such that $T$  contains  a subgraph isomorphic to  $H$.  The \dfn{planar anti-Ramsey number of $H$}, denoted  $ar_{_\mathcal{P}}(n, H)$,  is the  maximum number  of colors in an edge-coloring of  a plane triangulation $T\in \mathcal{T}_n(H)$ such that $T$ contains   no    rainbow copy of $H$.   Analogous to the relation between anti-Ramsey numbers and Tur\'an numbers proved in~\cite{ESS},  planar anti-Ramsey numbers are closely related to planar Tur\'an numbers~\cite{LSS2}, where  the \dfn{planar Tur\'an number of $H$}, denoted  $ex_{_\mathcal{P}}(n,H)$, is  the maximum number of edges of  a planar graph on $n$ vertices that contains no subgraph isomorphic to $H$.  
 
 \begin{prop}[\cite{LSS2}]\label{LU} Given a planar graph $H$ and a positive integer $n\ge |H|$, 
 $$1+ex_{_\mathcal{P}}(n,\mathcal{H})\le ar_{_\mathcal{P}}(n,H) \le ex_{_\mathcal{P}}(n, H),$$
 where $\mathcal{H}=\{H-e: \, e\in E(H)\}$. 
 \end{prop}
 
    Dowden~\cite{Dowden} began the study of planar Tur\'an numbers  (under the name of ``extremal" planar graphs). Results on planar Tur\'an numbers of paths and cycles can be found in~\cite{Dowden, LSS1}. The  study of planar anti-Ramsey numbers   was initiated by Hor\v{n}\'ak,  Jendrol$'$,  Schiermeyer and  Sot\'ak~\cite{HJSS} (under the name of rainbow numbers). Results on planar anti-Ramsey numbers of paths and cycles can be found in~\cite{HJSS, LSS2}. 
Colorings of plane graphs  that avoid rainbow faces have also been  studied,  see, e.g., \cite{Kral, JMSS, West, Zykov}. Various results on anti-Ramsey numbers can be found in:  \cite{Alon, Jiang2004,   Jiang2002, Jiang2009, Li, M, Ingo} to name a few. \medskip

Finding  exact values of $ar_{_\mathcal{P}}(n, H)$ is far from trivial. As observed in~\cite{HJSS}, an induction argument in general  cannot be applied to compute $ar_{_\mathcal{P}}(n, H)$ because  deleting a vertex from a plane triangulation may result in a graph that is no longer a plane triangulation.  
In this paper, we study  planar anti-Ramsey numbers of matchings. For all $t\ge1$,  let $M_t$ denote a matching of size $t$.  In \cite{M4}, the exact value   of $ar_{_{\mathcal{P}}}(n,    {M}_t)$ when $t\le4$  was  determined, and   lower and  upper bounds for  $ar_{_{\mathcal{P}}}(n,    {M}_t)$ were also established for all $t\ge5$ and $n\ge2t$.
Recently, the exact value  of $ar_{_{\mathcal{P}}}(n,    {M}_5)$ was determined  in \cite{M5} and an improved upper bound for  $ar_{_{\mathcal{P}}}(n,    {M}_t)$  was also  obtained in \cite{M5}. We summarize the results in \cite{M4, M5} below.

\begin{thm}[\cite{M4}]\label{M4}  Let $n$ and $t$ be positive integers.  Then
  \begin{enumerate}[(a)]

   \item for all    $n\ge 7$,  $ar_{_{\mathcal{P}}}(n,    {M}_3)= n $.
    \item for all    $n\ge 8$,  $ar_{_{\mathcal{P}}}(n,    {M}_4)= 2n-2$.
    \item  for all    $t\ge5$ and $n\ge 2t$, $ 2n+2t -10 \le ar_{_{\mathcal{P}}}(n,    {M}_t) \le 2n+2k-7+2{{2t-2}\choose 3}$.
     \end{enumerate}
\end{thm}

\begin{thm}[\cite{M5}]\label{M5} Let $n$ and $t$ be positive integers.  Then
  \begin{enumerate}[(a)]
    \item for all    $n\ge 11$,  $ar_{_{\mathcal{P}}}(n,    {M}_5)= 2n $.
    \item for all $t\ge5$ and $n\ge2t$, $ar_{_{\mathcal{P}}}(n,M_t)\le 2n+6t-17$.
    \end{enumerate}
\end{thm}

In this paper,  we further improve the existing lower   and upper bounds   for  $ar_{_{\mathcal{P}}}(n,    {M}_t)$.

\begin{thm}\label{LUB}
For all $t\ge6$  and   $n\geq 3t-6$,    $2n+3t-15\le ar_{_{\mathcal{P}}}(n,    {M}_t)\le  2n+4t-14$.
\end{thm}

  Theorem~\ref{LUB}   significantly improves the lower bound   in Theorem~\ref{M4}(c)  and the new upper bound in Theorem~\ref{M5}(b).   We believe that for each $t\ge6$, the lower bound we obtained in Theorem~\ref{LUB} is the exact value  of $ar_{_{\mathcal{P}}}(n,    {M}_t)$ for sufficiently large $n$.   This is indeed the case for $M_6$.  
 
\begin{thm}\label{M6}
  For all $n\ge 30$,    $ar_{_{\mathcal{P}}}(n,    {M}_6)= 2n+3$.
\end{thm}

   It seems that the method we developed  in the proof of Theorem~\ref{M6} can be applied to close the gap in Theorem~\ref{LUB}.  We prove Theorem~\ref{LUB} in Section~\ref{main}  and Theorem~\ref{M6} in Section~\ref{main*}.

\section{Proof of Theorem~\ref{LUB}}\label{main}

\begin{figure}[htbp]
\centering
\includegraphics*[scale=0.3 ]{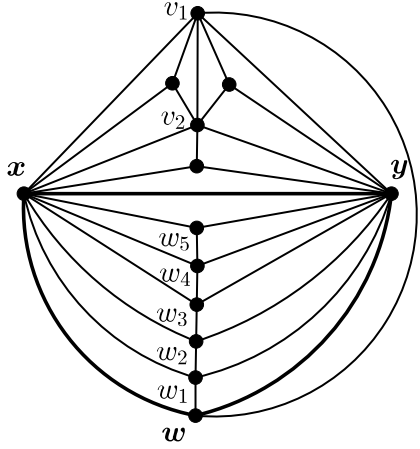}
\caption{The construction of $T$ when $t=6$ and $n=13$}\label{T}
\end{figure}

We are ready to  prove Theorem~\ref{LUB}.     Let $t,n$ be given as in the statement. We first prove that  $ar_{_{\mathcal{P}}}(n,    {M}_t)\ge 2n+3t-15$.   Let $P$ be a path with vertices $v_1, v_2, \dots, v_{t-4}$ in order. Let $H$ be the plane triangulation obtained from $P$ by adding two adjacent vertices $x, y$ and joining each of $x$ and $y$ to all the vertices on $P$ with the outer face of $H$ having  vertices $x,y,v_1$ on its boundary. Then $|H|=t-2\ge4$ and $H$ is hamiltonian. Let $T_{_H}$ be the plane triangulation  obtained from $H$ by adding a new vertex to each  face $F$ of $H$ and   then joining it to all vertices on the boundary of $F$.  Then $T_{_H}$ is a plane triangulation on $(t-2)+(2(t-2)-4)=3t-10$ vertices. Let $w$ be the new vertex added to the outer-face of $H$.  Let $T$ be the plane triangulation on $n$ vertices obtained from $T_{_H}$ by  adding $n-(3t-10)\ge4$ vertices, say $w_1, w_2, \ldots, w_{_{n-3t+10}}$, to the face   of $T_H$ containing $x,y, w$,  such that $ww_1, w_{_{n-3t+10}}x, w_{_{n-3t+10}}y\in E(T)$, and for all $i\in[n-3t+9]$, $w_i$ is adjacent to $x,  y, w_{i+1}$ in $G$.  The construction of $T$ when $t=6$ and $n=13$ is depicted in Figure~\ref{T}.  Clearly, $T\in \mathcal{T}_n(M_t)$. Let $c$ be an edge-coloring of $T$ by first coloring all the edges  $ww_1, w_1w_2, \ldots, w_{_{n-3t+9}} w_{_{n-3t+10}}$  by color $1$ and then all the remaining edges of $T$ by distinct colors other than $1$. It can be easily checked that $T$ has no rainbow $M_t$ under the coloring $c$ and the total number of colors   used  by $c$ is $(3n-6)-(n-3t+10)+1=2n+3t-15$. This proves that   $ar_{_{\mathcal{P}}}(n,    {M}_t)\ge 2n+3t-15$, as desired. \medskip

It remains to prove that    $ar_{_{\mathcal{P}}}(n, M_t)\le 2n+4t-14$.
 Suppose $ar_{_{\mathcal{P}}}(n, M_t)\ge 2n+4t-13$ for some $t\ge6$ and  $n\geq 3t-6$.    Then  there exists a $T \in\mathcal{T}_n(M_t)$  such that $T$ has no rainbow copy of $M_t$ under some onto mapping     $c:E(T)\rightarrow [k]$, where $k\ge 2n+4t-13$. We choose such a $T$ with $t$ minimum.      Let   $G$ be a rainbow spanning subgraph of $T$ with $k$ edges. Then $G$ does not contain  $M_t$ because $T $ has no rainbow copy of $M_t$.  By minimality of $t$ and Theorem~\ref{M5}(a) (when $t=6$),   $G$ contains a copy of $M_{t-1}$.
Let $M:=\{u_iw_i\in E(G):    i\in[ t-1]\}$ be a matching of size $t-1$ in $G$, and let $V(M):=\{u_1, \ldots, u_{t-1}, w_1, \ldots, w_{t-1}\}$.     Let $R:=V(G)\less V(M)$.  For each $i\in[t-1]$, we may assume that $|N_G(u_i)\cap R|\le |N_G(w_i)\cap R|$.   Since $M$ is the largest matching in $G$, we see that   $G$ has no   $M$-augmenting path. It follows that    $G[R]$ has no edges,  and      for each $i\in[t-1]$,  
either  $|N_G(u_i)\cap R|=0$ or $|N_G(u_i)\cap R|=1$ with $N_G(u_i)\cap R =N_G(w_i)\cap R$.  We may further assume that $u_1w_1, \ldots, u_\ell w_\ell$ are such that $ |N_G(u_i)\cap R|=1 $  for all $0\le i\le \ell$, and $|N_G(u_j)\cap R|=0 $ for all $j\in\{\ell+1, \ldots, t-1\}$, where $0\le \ell\le t-1$. Then $e_G(\{u_1, \ldots, u_\ell, w_1, \ldots, w_\ell\}, R)=2\ell$. Let $L:=\{w_{\ell+1}, \ldots, w_{t-1}\}$.   Then   $e_G(L, R)=0\le2n-2t-2\ell-2$ when $\ell=t-1$, and  $e_G(L, R)\le 2(n-(t-1)-\ell)-4=2n-2t-2\ell-2$ when $\ell\le t-2$ because  $G[L\cup  R]\less E(G[L])$ is a planar bipartite  graph on $n-(t-1)-\ell \ge  3$ vertices.  Since $G$ is planar and  $|V(M)|=2(t-1)>3$, we have  $e(G[V(M)])\le 3(2t-2)-6=6t-12$. Thus, $ e(G)= e(G[V(M)])+e_G(\{u_1, \ldots, u_\ell, w_1, \ldots, w_\ell\}, R)+e_G(L, R)\le (6t-12)  +2\ell+(2n-2t-2\ell-2)=2n+4t-14$, contrary to $e(G)=k\ge 2n+4t-13$. This completes the proof of Theorem~\ref{LUB}. \qed\\

\noindent {\bf Remark.} For $t\ge7$, the condition ``$n\ge3t-6$'' in the statement of Theorem~\ref{LUB} can be replaced by ``$n\ge3t-7$''. \medskip

\section{Proof of Theorem~\ref{M6}}\label{main*}

We need to introduce more notation that shall be used in this section only. 
For   $n\ge 3$,  let  $\mathcal{T}_n$ be the set of all plane triangulations   on $n$ vertices, and let  $\mathcal{T}^-_n$  be the set of all planar graphs with $n$ vertices and $3n-7$ edges.  Clearly, every graph in $\mathcal{T}^-_n$  is isomorphic to a plane triangulation on $n$ vertices with one edge removed.  By abusing notation, let  $e(\mathcal{T}_n): =3n-6$  and $e(\mathcal{T}^-_n): =3n-7$.  
 It is known that    every plane triangulation  on $n\ge4$ vertices  is $3$-connected. It is also known that every plane triangulation  on $  n\le10$ vertices  has a Hamilton cycle   and every plane triangulation  on $  n\ge11$ vertices does not   necessarily have a Hamilton cycle\footnote{The third author would like to thank Jason Bentley, a Ph.D. student at the University of Central Florida, for his help in carefully verifying these facts with her.}. We summarize these facts as follows.

\begin{ob}\label{triangulations} Let $T$ be a   planar triangulation on   $n\ge4$ vertices. Then
\begin{enumerate}[(a)]
\item $T$   is $3$-connected.
\item for every $n\le10$,   $T$ has a Hamilton cycle.
\item  for every $n\ge11$,   $T$ does not necessarily have a Hamilton cycle.
\end{enumerate}
\end{ob}

Let $o(H)$ denote  the  number of  odd   components in a graph $H$.  We shall make use of the following theorem in the proof of Theorem~\ref{M6}.

\begin{thm}[Berge-Tutte Formula]\label{Tutte}
    Let $G$ be a graph on $n$ vertices and  let $d$ be the size of a maximum matching of $G$. Then there exists  an $S\subseteq V(G)$ with    $|S|\le d$   such that $$2d= n-o(G\less S)+|S|.$$  
    Moreover,  each odd   component of $G\less S$ is factor-critical.
\end{thm}
\bigskip

\noindent{\bf Proof of Theorem~\ref{M6}}:   Let $n\ge30$ be an integer. By Theorem~\ref{LUB}, $ar_{_{\mathcal{P}}}(n,    {M}_6)\ge 2n+3$. We next show that $ar_{_{\mathcal{P}}}(n,    {M}_6)\le 2n+3$. Suppose $ar_{_{\mathcal{P}}}(n,    {M}_6)\ge 2n+4$. Then there exists a   $T\in\mathcal{T}_n(M_6)$   such that $T$ has no rainbow $M_6$ under some   onto mapping $c:E(T)\rightarrow [k]$, where $k\ge 2n+4$.  Let $G$ be a  rainbow spanning subgraph of $T$ with $k$ edges.
 By Theorem~\ref{M5}(a),   $G$ has a  copy of $M_5$.   Clearly, $G$ has no copy of $M_6$ because $T$ has no rainbow copy of $M_6$ under $c$. By Theorem~\ref{Tutte}, there exists an  $S\subseteq V(G)$ with $s:=|S|\le 5$ such that
   $q:=o(G\setminus S)=n+s-10$. 
    Let $H_1,H_2,\ldots,H_q$ be all the odd components of $G\setminus S$.   We may assume that $|H_1|\le |H_2|\le\cdots\le|H_q|$. Let $r: =\max\{i:\, |H_i|=1\}$.  Then $n=|G|\ge |S|+(| H_1|+\cdots+ |H_r|) +(| H_{r+1}|+ \cdots   |H_q|)\ge s+r+3(q-r)=s-2r+3(n+s-10)$. It follows that  $r\ge  n+2s-15\ge15$.   Let $S:=\{v_1,\ldots,v_{s}\}$ when $s\ge1$ and $ V(H_i)=\{u_i\}$ for all $i\in[r]$.  We may further assume that $d_{G}(u_1)\ge d_{G}(u_2)\ge \cdots\ge d_{G}(u_r)$.   Let $U : =\{u_1, \ldots, u_r \}$ and   $W: = V(G)\setminus (S\cup U)$.  Then $w:=|W|=n-s-r$ and  $e_G(U, S)\le 2(r+s)-4$ when $r+s\ge3$.  We next  prove several claims.\\

\noindent {\bf Claim 1.}\,   If $G$  has   two edge-disjoint matchings of size $5$, say   $M'$ and $M''$,   then  $T[V(G)\less V(M'\cup M'')] $ has no edges.

\pf Suppose $T[V(G)\less V(M'\cup M'')]$ has an edge $e$. We may assume that  $c(e)\ne c(e')$ for all $e'\in M'$. But then $M'\cup\{e\}$ is a rainbow $M_6$ in $T$ under the coloring $c$, a contradiction. \qed\\

\noindent {\bf Claim 2.}\, If $w+s\le9$, then $H:=G[ W\cup S\cup\{u_1, \ldots, u_{_{10-w-s}}\}]\notin \mathcal{T}_{10}$.

\pf Suppose  $H\in \mathcal{T}_{10}$. Then     $H$ has a Hamilton cycle by Observation~\ref{triangulations}(b), and thus  has two edge-disjoint  matchings of size $5$. By Claim 1, $T[\{u_{_{11-w-s}}, \ldots, u_r\}]$ has no edges. But then
\begin{align*}
 e(T )&=e(T[ W\cup S\cup  \{u_1,\ldots,u_{_{10-w-s}}\}])+e_T(\{u_{_{11-w-s}}, \ldots, u_r\}, W \cup S\cup\{u_1, \ldots, u_{_{10-w-s}}\})\\
& \le e(\mathcal{T}_{10})+(2n-4)=24+(2n-4),
\end{align*}
  which implies that  $n\le 26$ because $e(T)= 3n-6$, contrary to  $n\ge30$. \qed\\

\noindent {\bf Claim 3.}\, $|H_q|\ge3$.

\pf Suppose $|H_q|<3$. Then $r=q$ and   so  $w+s=n-q  =10- s$.  It follows that 
\begin{align*}
 2n+4\le e(G)&=e(G[W\cup S])+e_G(U, S)\le (3(w+s)-6)+(2(n-w)-4)\\
 &=2n+w+3s-10=2n+(10-2s)+3s-10=2n+s,
\end{align*}
 which implies that $s\ge4$. If  $s=4$, then    $w=2$. But then 
 \[2n+4\le e(G)= e(G[S])+e(G[W])+e_G(S, U\cup W)\le e(\mathcal{T}_{4})+1+2n-4=2n+3,\]
  which is impossible.
Thus   $s=5$, and so  $w+s=5$ and   $r=q=n-5$.  By Claim 2,  $e(H)\le 23$.  
Since $2n+4\le e(G)=e(G[S])+e_G(U, S)\le e(G[S])+2n-4$, we see that $e(G[S])\ge8$.  Thus $G[S]\in \mathcal{T}_5$ or $G[S]\in \mathcal{T}_5^{-}$.
  Then  $d_{G}(u_1)\le 12-e(G[S])$, else $e(G[S\cup\{u_1\}])\ge 13>e(\mathcal{T}_6)=12$, a contradiction.
Suppose  $G[S]\in \mathcal{T}_5$.      Then $d_{G}(u_1)\le 3$ and $e_G(U, S)=2n-5$.  It follows that   $d_G(u_5)=3$, else  $e_G(U, S)=e_G(\{u_1, u_2, u_3, u_4\}, S)+e_G(\{u_5, \ldots, u_r\}, S)\le 12+2(n-9)=2n-6$.  But  then  $e(H)  =e(G[S])+e_G(\{u_1, \ldots, u_5\}, S)=9+3\times 5=24$, contrary to $e(H)\le23$.
  This proves that  $G[S]\in \mathcal{T}_5^{-}$.  We may assume that $v_1v_2\notin E(G[S])$.  Then   $d_{G}(u_1)\le 12-e(G[S])=4$, and  $e_G(U, S)=e(G)-e(G[S])\ge (2n+4)-8=2n-4$, which implies that  $e_G(U, S)=2n-4$. 
    If  $d_{G}(u_2)\ge4$, then $e(G[S\cup \{u_1,u_2\}])=8+8=16>e(\mathcal{T}_7)$, a contradiction. Thus $d_{G}(u_2)\le 3$. Furthermore, $d_G(u_5)=3$, else  $e_G(U, S)=e_G(\{u_1, u_2, u_3, u_4\}, S)+e_G(\{u_5, \ldots, u_r\}, S)\le (4+9)+2(n-9)=2n-5$, contrary to $e_G(U, S)=2n-4$.    Since $e(H)\le 23$, we see that  $d_{G}(u_1)=3$. Then $d_G(u_6)=3$, else $e_G(U, S)\le 15+2(n-10)=2n-5$, a contradiction. Since $G[S\cup\{u_1, \ldots, u_6\}]$ does not contain $K_{3,3}$ as a subgraph, we may assume that $v_1u_6\in E(G)$. Then
  $H\in \mathcal{T}_{10}^-$ because $e(H)=e(G[S])+e_G(\{u_1, \ldots, u_5\},S)=8+15=e(\mathcal{T}_{10}^-)$. Note that  $H+v_1v_2\in \mathcal{T}_{10}$. By Observation~\ref{triangulations}(b),   $H$ has a hamiltonian  path with   $v_1 $ as an end. Since $v_1u_6\in E(G)$, we see that   $G[S\cup\{ u_1, \ldots, u_6\}]$ has two edge-disjoint   matchings of size $5$. By Claim 1, $T[\{u_7, \ldots, u_r\}]$ has no edges. But then
\begin{align*}
3n-6=e(T )&=e(T[S\cup \{u_1,\ldots,u_6\}])+e_T(\{u_7, \ldots, u_r\}, S\cup \{u_1, \ldots,u_6\})\\
 &\le e(\mathcal{T}_{11})+(2n-4)=27+(2n-4),
 \end{align*}
  which implies that  $n\le 29$, contrary to  $n\ge30$. \qed\\

  By Claim 3, $w\ge3$ and $r\le q-1$. Then $ n\ge s+|H_1|+\cdots+|H_q|\ge s+q+2=n+2s-8$, which implies that $s   \le4$, with $s=4$ only when $w=|H_q|=3$.   \\

  \noindent {\bf Claim 4.}\, $|H_{q-1}|\ge3$.

  \pf Suppose  $|H_{q-1}|=1$. By   Claim 3,   $r= q-1$.  Thus     $w +s=n-r =n-(q-1)   =11- s$.     It follows that 
  \begin{align*}
  2n+4\le e(G)&=e(G[W\cup S])+e_G(U, S)\le (3(w+s)-6)+(2(n-w)-4)\\
  &=2n+w+3s-10=2n+(11-2s)+3s-10=2n+s+1,
  \end{align*}
   which implies that $s\ge3$.   If  $s=3$, then    $w+s=8$ and so $d_G(u_1) =3$, else   $e(G)\le e(G[W\cup S])+e_G(U,S)\le e(\mathcal{T}_8) + 2(n-8)=2n+2$, a contradiction. Since  $G$ does not contain $K_{3,3}$ as a subgraph, we see that     $d_G(u_3)\le2$.  By Claim 2,  $e(H)\le23$. Thus $ e(G[W\cup S])= e(H)-d_G(u_1)-d_G(u_2)\le 20- d_G(u_2) $.  Note that  $e_G(U, S)\le 3+d_G(u_2)+2(n-w-s-2)=2n-17+d_G(u_2)$.  But then 
   \[2n+4\le e(G)= e(G[W\cup S]) +e_G(U, S)\le (20-d_G(u_2))+(2n-17+d_G(u_2))=2n+3,\]
    which is impossible.  Thus $s=4$. Then $w=3$. It follows that    $W=V(H_q)$ and $G[W]=K_3$ because  $G[W]$ is factor-critical.   Then $d_G(u_2)\ge3$, else  $ e_G(U, S)\le 4+2(n-8)=2n-12$ and so $2n+4\le e(G)= e(G[W\cup S]) +e_G(U, S)\le e(\mathcal{T}_7)+(2n-12)=2n+3$, a contradiction. Since $e_G(U, S)\le 2(n-3)-4=2n-10$, we see that  $  e(G[W\cup S])=e(G)- e_G(U, S)\ge (2n+4)- (2n-10) \ge14$. Thus $G[W\cup S]\in \mathcal{T}_7$ or $G[W\cup S]\in \mathcal{T}_7^-$.
  Suppose $G[W\cup S]\in\mathcal{T}_7$. Then $d_G(u_i)\le3$ for all $i\in [r]$. Since   $e_G(U, S)=e(G)  - e(\mathcal{T}_7)\ge 2(n-7)+3$, we see that $d_G(u_1) =d_G(u_2) =d_G(u_3) =3$. But then $e(H)=24$,  contrary to $e(H)\le23$. This proves that   $G[W\cup S]\in\mathcal{T}_7^-$. Then $d_G(u_2)\le3$, else  $e(G[W\cup S\cup\{u_1, u_2\}])=e(\mathcal{T}_7^-)+8\ge 14+8>e(\mathcal{T}_9)$, a contradiction.  Since   $e_G(U, S)=e(G)  - e(\mathcal{T}^-_7)\ge 2(n-7)+4$, we see that $d_G(u_2) =d_G(u_3) =3$. Then 
     $d_G(u_1) =3$    because $e(H)\le 23$.  Since $G[W\cup S]\in\mathcal{T}_7^-$ and $G[W ]=K_3$, we may assume that $G[W\cup S]+v_1v\in\mathcal{T}_7$ for some $v\in W\cup S$ with $v\ne v_1$.  Note that 
  $G[S\cup \{u_2, u_3, u_4\}]$ does not contain $K_{3,3}$ as a subgraph. We may further assume that $v_1u_4\in E(G)$. Then  $H+v_1v\in\mathcal{T}_{10}$. By Observation~\ref{triangulations}(b),  $H$ has a hamiltonian path with $v_1$ as an end. Since $v_1u_4\in E(G)$, we see that  $G[W\cup S\cup \{u_1, u_2, u_3, u_4\}]$ has two edge-disjoint matchings of size $5$. By Claim 1,
  $T[\{u_5, \ldots, u_r\}]$ has no edges. But then
\begin{align*}
3n-6=e(T )&=e(T[W\cup S\cup \{u_1,u_2, u_3,u_4\}])+e_T(\{u_5, \ldots, u_r\}, W\cup S\cup \{u_1, u_2, u_3,u_4\})\\
 &\le e(\mathcal{T}_{11})+(2n-4)=27+(2n-4),
 \end{align*}
  which implies that  $n\le 29$, contrary to  $n\ge30$. \qed\\

  By Claim 4, $w\ge6$ and $r\le q-2$. Then $ n\ge s+|H_1|+\cdots+|H_q|\ge s+q+4=n+2s-6$, which implies that $s\le  3$,  with $s=3$ only when $w=|H_q|+|H_{q-1}|=3+3=6$ and $r=q-2$.   Suppose $s=3$. Then $w=6$.  By Claim 2, $e(H)\le23$.  Since  $G[S\cup\{u_1, u_2, u_3\}]$  does not contain $K_{3,3}$ as a subgraph, we see that    $d_G(u_3)\le2$.   Then $d_G(u_2)=3$ and $H\in \mathcal{T}_{10}^-$,  else  either
\[e(G )=e(H)+e_G(\{u_2, \ldots, u_r\},   S )\le 23+2(n-10)=2n+3,  \text{ or }\]  
\[e(G )=e(H)+e_G(\{u_2, \ldots, u_r\},   S )\le 22+3+2(n-11)=2n+3,\] contrary to $e(G)\ge 2n+4$ in both cases.  Let $V(H_{q-1})=\{x, y, z\}$ and  $V(H_{q})=\{x', y', z'\}$. Since  $H\in\mathcal{T}_{10}^-$, we may assume that $xv_2, x'v_3\in E(G)$.  It follows that $G[W\cup S\cup\{u_1, u_2\}]$ has two edge-disjoint matchings of size $5$, namely, $\{v_1u_1,  v_2x, v_3u_2, yz, x'y'\}$ and $\{v_1u_2, v_2u_1, v_3x', xy, y'z'\}$.  By Claim 1,  $T[\{u_3, \ldots, u_r\}]$ has no edges. But then
\begin{align*}
3n-6=e(T )&=e(T[W\cup S\cup \{u_1,u_2  \}])+e_T(\{u_3, \ldots, u_r\}, W\cup S\cup \{u_1, u_2  \})\\
 &\le e(\mathcal{T}_{11})+(2n-4)=27+(2n-4),
 \end{align*}
 which implies that  $n\le 29$, contrary to  $n\ge30$.  Thus  $s\le 2$. Then  $d_G(u_i)\le s$ for all $i\in[r]$ and    $S$  is a vertex-cut of  $G[S\cup W]$   because $r\le q-2$.  Then $e_G(U, S)\le rs$ and  $e(G[S\cup W])\le e(\mathcal{T}_{w+s})-2=3(w+s)-8=3(n-r)-8$. But then
\begin{align*}
2n+4\le e(G)&=e_G(U, S)+ e(G[S\cup W])\\
&\le rs +3(n-r)-8\\
&=3n-(3-s)r-8\\
&\le 3n-(3-s)(n+2s-15)-8,
\end{align*}
which is impossible because    $s\le 2$ and     $n\ge30$.   This completes the proof of Theorem~\ref{M6}.   \qed\\

 \noindent {\bf Remark.} In the proof of Theorem~\ref{M6},  Claim 1 is applied to two vertex-disjoint matchings, instead of edge-disjoint matchings. It seems that the method we developed in the proof of Theorem~\ref{M6} can be used to close the gap in Theorem~\ref{LUB}.  \\


\noindent {\bf Acknowledgments.}  The authors would like to thank Jingmei Zhang for helpful comments. \\
 Gang Chen  would like to thank the University of Central Florida for hosting his visit. His research  is partially supported by NSFC under the grant number 71561022 and Overseas Training Program for Faculty at   Ningxia University.

\end{document}